\documentclass[a4paper,reqno,10pt,oneside]{amsart}
\usepackage{amsmath}
\usepackage[normalem]{ulem}
\usepackage[left=2.61cm,right=2.61cm,top=2.72cm,bottom=2.72cm]{geometry}
\usepackage[colorlinks=true,urlcolor=blue,
citecolor=red,linkcolor=blue,linktocpage,pdfpagelabels,
bookmarksnumbered,bookmarksopen]{hyperref}
\usepackage[english]{babel}
\usepackage{graphicx}
\usepackage[toc,page,header]{appendix}
\usepackage{mathrsfs}
\usepackage{amssymb,amsmath, amsfonts, amsthm}
\usepackage{latexsym}
\usepackage{enumitem}

\allowdisplaybreaks[1]

\numberwithin{equation}{section}

\newtheorem{theorem}{Theorem}[section]
\newtheorem{proposition}[theorem]{Proposition}
\newtheorem{corollary}[theorem]{Corollary}
\newtheorem{lemma}[theorem]{Lemma}
\theoremstyle{definition}
\newtheorem{remark}[theorem]{Remark}
\theoremstyle{definition}
\newtheorem{definition}[theorem]{Definition}
\theoremstyle{definition}

\renewcommand{\le}{\leqslant}
\renewcommand{\ge}{\geqslant}

\newcommand{\N}{\mathbb{N}}
\renewcommand{\S}{\mathbb{S}}

\newcommand{\beq}{\begin{equation}}
\newcommand{\eeq}{\end{equation}}
\newcommand{\beqs}{\begin{equation*}}
\newcommand{\eeqs}{\end{equation*}}
\newcommand{\beqn}{\begin{eqnarray}}
\newcommand{\eeqn}{\end{eqnarray}}
\newcommand{\beqns}{\begin{eqnarray*}}
\newcommand{\eeqns}{\end{eqnarray*}}
\newcommand{\bdoc}{\begin{document}}
\newcommand{\edoc}{\end{document}}
\newcommand{\be}{\begin{enumerate}}
\newcommand{\ee}{\end{enumerate}}
\newcommand{\bdescr}{\begin{description}}
\newcommand{\edescr}{\end{description}}
\newcommand{\ba}{\begin{array}}
\newcommand{\ea}{\end{array}}
\newcommand{\intR}{\int_{\mathbb R^N}}
\newcommand{\R}{\mathbb R}
\newcommand{\B}{\mathbb B}
\newcommand{\C}{\mathbb C}
\renewcommand{\H}{\mathcal H}
\renewcommand{\L}{\mathbb L}
\newcommand{\parallelsum}{\mathbin{\!/\mkern-5mu/\!}}
\newcommand{\e}{\varepsilon}
\newcommand{\SD}{\Sigma_D}

\renewcommand{\appendixpagename}{\centering Appendix}

\newcommand{\todo}[1]{\text{\colorbox{yellow}{#1}}}

\newcommand{\modified}[1]{\text{\colorbox{green}{#1}}}

\begin{document}

\title[Bifurcation from bubbles in nonconvex cones]{Bifurcation from bubbles in nonconvex cones}

\date{\today}

\author{Filomena Pacella}
\address{F. Pacella. Dipartimento di Matematica ``Guido Castelnuovo'', Sapienza Universit\`a di Roma,  P.le Aldo Moro 2, 00185 Roma, Italy}
\email{filomena.pacella@uniroma1.it}

\author{Camilla Chiara Polvara}
\address{C.C. Polvara. Dipartimento di Matematica ``Guido Castelnuovo'', Sapienza Universit\`a di Roma,  P.le Aldo Moro 2, 00185 Roma, Italy}
\email{camilla.polvara@uniroma1.it}

\author{Luigi Provenzano}
\address{L. Provenzano.  Dipartimento di Scienze di Base e Applicate per l'Ingegneria ``SBAI'', Sapienza Universit\`a di Roma, Via Antonio Scarpa, 14, 00161 Roma, Italy}
\email{luigi.provenzano@uniroma1.it}

\subjclass{35B32; 35B33; 35J60}
\keywords{Critical exponent equation, Neumann problem in cones, Bifurcation, Symmetry breaking}
\thanks{\emph{Acknowledgements.} Research partially supported by Gruppo Nazionale per l'Analisi Matematica, la Pro\-ba\-bi\-li\-t\`a e le loro Applicazioni (GNAMPA) of the Istituto Nazionale di Alta Matematica (INdAM). F.P. and C.C.P. have been supported by PRIN 2022 PNRR project 2022AKNSE4-Next Generation EU “Variational and Analytical aspects of Geometric PDEs”, founded by the European Union. L.P. has been supported by the project ``Perturbation problems and asymptotics for elliptic differential equations: variational and potential theoretic methods'' funded by the European Union – Next Generation EU and by MUR-PRIN-2022SENJZ3 and by the project ``Analisi Geometrica e Teoria  Spettrale su varietà Riemanniane ed Hermitiane'' of the INdAM GNSAGA}

\begin{abstract}
 We investigate the
 Neumann problem for the critical semilinear elliptic equation in cones. The standard bubble provides a family of radial solutions, which are known to be the only positive solutions in convex cones. For nonconvex cones, symmetry breaking may occur and the symmetry breaking is related to the first nonzero Neumann eigenvalue of the Laplace–Beltrami operator on the domain $D\subset\S^{N-1}$, that spans the cone.
We construct a one-parameter family of domains on the sphere whose first eigenvalue crosses the threshold at which the bubble loses stability. Under the assumption that this eigenvalue is simple, we prove, via the Crandall–Rabinowitz bifurcation theorem, the existence of a branch of nonradial solutions bifurcating from the standard bubble. Moreover we show that the bifurcation is global.
\end{abstract}

\maketitle

\section{Introduction}
Let $D$ be a smooth domain on the unit sphere $\S^{N-1}\subset\mathbb R^N$, with $N\ge 3,$ and let $\Sigma_D\subset\R^N$ be the cone spanned by $D$:
\beq\label{cono}
\Sigma_D=\{x\in\R^N: x=\rho p,p\in D,\rho\in(0,+\infty)\},
\eeq
and assume that $\Sigma_D$ is a Lipschitz domain.
We consider the following semilinear critical Neumann problem
\beq\label{printro}
\begin{cases}
    -\Delta u=u^{2^*-1}&\text{in }\Sigma_D,\\
    \frac{\partial u}{\partial\nu}=0&\text{on }\partial\Sigma_D\setminus\{0\},\\
    u>0&\text{in }\Sigma_D,
\end{cases}
\eeq
where $2^*=\frac{2N}{N-2}$ is the critical Sobolev exponent for the immersion $\mathcal{D}^{1,2}(\Sigma_D)$ into $L^{2^*}(\Sigma_D)$ and
$$
\mathcal{D}^{1,2}(\Sigma_D)=\{u\in L^{2^*}(\Sigma_D):|\nabla u|\in L^2(\Sigma_D)\}.
$$
It is well known that the radial function:
\beq\label{bubbleintro}
U(x)=c_0\bigg(\frac{1}{1+|x|^2}\bigg)^\frac{N-2}{N},\qquad x\in\Sigma_D,
\eeq
is a solution of \eqref{printro}, for $c_0=(N(N-2))^\frac{N-2}{4}$, as well as any of its rescaling. We recall that $U$ is a minimizer for the Sobolev quotient in the whole $\R^N$ and we refer to it as the standard bubble. Moreover, in the case of $\R^N$, the function $U$, its rescalings and translations are the only positive solutions of \eqref{printro}. Obviously they are also the only radial solutions of \eqref{printro} and we refer to them as the trivial solutions of \eqref{printro}. Then the question is whether in $\Sigma_D$, for $D\ne\S^{N-1}$, there exist also nonradial solutions. 

The study of solutions of the critical equation in cones dates back to 1988 where it arised in connection with the Sobolev and Isoperimetric inequality in cones (see \cite{LP,LPT}).
In particular in \cite{LPT} it was proved that if the cone $\Sigma_D$ is convex then the bubbles are the only positive solution of \eqref{printro}. More recently this symmetry result was extended to critical equations for more general operators in \cite{CFR}. Thus the question is whether nonradial solutions of \eqref{printro} exist in nonconvex cones.

Since minimizers of the Sobolev quotient 
\beq\label{Sobqintro}
Q_D(u)=\frac{\big(\int_{\Sigma_D}|\nabla u|^2\big)^\frac{1}{2}}{\big(\int_{\Sigma_D} |u|^{2^*}\big)^\frac{1}{2^*}},\qquad u\in\mathcal{D}^{1,2}(\Sigma_D),u\ne 0,
\eeq
produce solutions of \eqref{printro}, a first result was obtained in \cite{CP} by showing that nonradial minimizers exist in some nonconvex cones. Later a more precise characterization of a class of nonconvex cones for which \eqref{printro} admits a nonradial solution (as a minimizer of \eqref{Sobqintro}) was obtained in \cite{CPP}.
It pointed out the crucial role played, in the break of symmetry, by the eigenvalue $\lambda_1(D)$ which is the first nonzero eigenvalue of the Laplace-Beltrami operator $-\Delta_{\S^{N-1}}$ on the domain $D$, with homogeneous Neumann boundary conditions.
The result of \cite{CPP} is the following:
\begin{theorem}\label{thvecchio}\cite{CPP}
    Let $\Sigma_D$ be a cone such that $\bar D\subset \S^{N-1}_+$ and $\lambda_1(D)<N-1$, where $\S^{N-1}_+$ is the half-sphere. Then the minimizers of the Sobolev quotient are nonradial and hence there exist nonradial solutions $v$ of \eqref{printro} which are also fast decaying, i.e.
    $$
    v(x)=O(|x|^{2-N}) \quad \text{as }|x|\to\infty.
    $$
\end{theorem}
It is interesting to quote that the same bound on the eigenvalue $\lambda_1(D)$, i.e. $\lambda_1(D)<N-1$, is also the key point to show breaking of symmetry results for overdetermined and constant mean curvature problems in cones (see \cite{IPW}). 

The proof of Theorem \ref{thvecchio} relies on a careful analysis of the Morse index of the standard bubble $U$, which shows that it becomes an unstable critical point of \eqref{Sobqintro} as soon as $\lambda_1(D)$ crosses the value $N-1$ (see Theorem 4.3 in \cite{CPP}). Instead when $\lambda_1(D)>N-1$ the bubble $U$ is a stable critical point for $Q$, i.e. it is a strict local minimum. Hence the analysis performed in \cite{CPP} shows that, varying the domains $D\subset \S^{N-1}$ (and hence the cone $\Sigma_D$) the bubble $U$ is degenerate for any cone $\Sigma_{D_1}$ such that $\lambda_1(D_1)=N-1$. This suggests that a bifurcation from the standard bubble could appear when $\lambda_1(D_1)=N-1.$ In other words a branch of nonradial solutions of \eqref{printro} in domains $D$ close to $D_1$ could emanate from the bubble $U$ in $D_1$.

This is the aim of the present paper. To state precisely our result, let $D=D_1$ be a domain in $\S^{N-1}_+$ such that $\lambda_1(D_1)=N-1$.
We consider the family of domains $D_\alpha$ which are diffeomorphic to $D_1$ through the diffeomorphism $\Phi_\alpha$ defined in Section \ref{sec2} (see \eqref{Phidef}), for $\alpha\in (0,\alpha^*)$, where $\alpha^*>0$ is the maximum value of $\alpha$ such that $\Phi_\alpha(D_1)$ is contained in $\S^{N-1}_+$. By construction, we have that $\alpha^*>1$.
We stress that, for every $\alpha$ the standard bubble $U$ is a solution of \eqref{printro} in $\Sigma_{D_\alpha}$, so we have the trivial curve of solutions $\alpha
\to(\alpha,U)$, for every $\alpha\in(0,\alpha^*)$. Our bifurcation result is the following.
\begin{theorem}\label{mainth}
Let $\Sigma_{D_\alpha}$ be a family of cones as above, $\alpha\in(0,\alpha^*)$ and assume that $\lambda_1({D_1})=N-1$ is a simple eigenvalue. Then the point $(1,U)$ is a bifurcation point for the trivial curve $\alpha\mapsto(\alpha,U)$ of solutions to \eqref{printro}. More precisely there exist $\epsilon>0$ and a $C^1$ curve  $(-\epsilon,\epsilon)\ni s\to (\alpha_s,v_s)$,  such that $\alpha_0=1$, $v_0=U$, $D_{\alpha_s}\ne D_{\alpha_1}$ and
$$
\begin{cases}
   -\Delta v_s=v_s^{2^*-1} &\text{in }\Sigma_{D_{\alpha_s}}\\
    v_s\in \mathcal{{D}}^{1,2}(\Sigma_{D_{\alpha_s}}),
\end{cases}$$
and $v_s$ is positive and nonradial.
\end{theorem}
To prove our result we will use the classical Crandall-Rabinowitz bifurcation theorem (see Theorem \ref{CRT}). The main difficulty in proving Theorem \ref{mainth} is to find a good family of domains $D_\alpha$ on $\S^{N-1}$ for which the corresponding eigenvalues $\lambda_1(D_\alpha)$ can be analyzed. Roughly speaking, the goal is to produce a family of domains for which $\lambda_1(D_{\alpha})$ crosses the value $N-1$ as $\alpha$ varies in $(0,\alpha^*)$ in a strict monotonic way. This is not easy, since the behavior of Neumann eigenvalues under domain perturbation can be quite involved, already in the Euclidean case; moreover, Neumann eigenvalues do not enjoy domain monotonicity and futhermore, on the sphere, we don't have the notion of homothety (hence we don't have simple transformations which allow to track explicitly the eigenvalues as in $\mathbb R^N$). In Section \ref{sec2} we construct a family of diffeomorphism $\Phi_\alpha$ on the hemisphere $\S^{N-1}_+$, which are, essentially, dilations with respect to the north pole. Then for the corresponding domain $D_\alpha=\Phi_\alpha(D)$ ($D$ a fixed domain suitably chosen) we are able to study the asymptotic behavior of $\lambda_1(D_\alpha)$, as $\alpha\to 0$, and, what is more important, to compute the derivative $\frac{d}{d\alpha}\lambda_1(D_\alpha).$

The proof that this derivative is negative at $\alpha=1$ (that is, when $\lambda_1(D_1)=N-1$) is the key point to get the so-called transversality condition in the application of Crandall-Rabinowitz theorem.

We also provide an example of a family of domains $D_\alpha$ for which $\lambda_1(D_\alpha)$ is simple (in a neighborhood of $\alpha=1$). However we point out that, generically, for almost all domains on $\S^{N-1}$ all Neumann eigenvalues are simple (see \cite{Henri}).

The other steps of the proof of Theorem \ref{mainth} follow by studying problem \eqref{printro} in a suitable space where, in particular, all functions are invariant by the Kelvin transform. This is inspired by the study of other critical exponent problems (see \cite{GGT}).

Finally we show that the bifurcation is global and that the Rabinowitz alternative holds. We refer to Section \ref{sec3} for the precise statement.
\bigskip

\textbf{Organization of the paper.} The paper is organized as follows. Section \ref{sec2} is devoted to the geometrical question of finding domains on the sphere to which the Crandall-Rabinowitz bifurcation Theorem can be applied. In Section \ref{sec3} we prove the bifurcation results.

\section{Geometric preliminaries, and a class of diffeomorphisms on the sphere}\label{sec2}

In this section we will construct a family $\{D_{\alpha}\}_{\alpha\in(0,\alpha^*)}$ of smooth domains contained in a hemisphere for which the first (non trivial) Neumann eigenvalue behaves as follows:
\begin{itemize}
\item $\lambda_1(D_1)=N-1$;
\item $\lambda_1(D_{\alpha})>N-1$ in $(0,1)$ and $\lambda_1(D_{\alpha})<N-1$ in $(1,\alpha^*)$
\item $\lambda_1(D_{\alpha})$ is simple in a neighborhood of $\alpha=1$;
\item $\frac{d}{d\alpha}\lambda_1(D_{\alpha})|_{\alpha=1}<0$.
\end{itemize}
We will also introduce a few geometric preliminaries on the Laplacian on a cone $\Sigma_D$ spanned by a spherical domains $D$, and its behavior under domain perturbation.

\subsection{Neumann eigenvalues on spherical domains}\label{ss3.1}

Let $\mathbb S^{N-1}, N\ge 3$, be the $(N-1)$-dimensional (unit) round sphere embedded in $\mathbb R^N$. By $\mathbb S^{N-1}_+$ we denote the upper hemisphere, namely $\mathbb S^{N-1}\cap\{x_N>0\}$. By $\Delta_{\mathbb S^{N-1}}$ we denote the Laplacian on $\mathbb S^{N-1}$ with respect to the (standard) round metric.

Let $D\subset\mathbb S^{N-1}_+$ be a smooth domain. By $\lambda_j(D)$ with $j\in \N$, we denote the eigenvalues of the Laplacian $\Delta_{\mathbb S^{N-1}}$ on $D$ with Neumann boundary conditions, i.e. we consider the following problem:
\begin{equation}\label{pb_lambda_def1}
\begin{cases}
-\Delta_{\S^{N-1}}Y_j = \lambda_j(D) Y_j & \text{on } D\\
\partial_{\nu_D} Y_j = 0 & \text{on } \partial D \,.
\end{cases}
\end{equation}
It is well known that $-\Delta_{\S^{N-1}}$ has a compact, self-adjoint resolvent in $L^2(D)$ and admits a sequence of eigenvalues 
\begin{equation} \label{lambdabe}
0=\lambda_0(D)<\lambda_1(D)\le \cdots \le \lambda_j(D) \le \cdots\nearrow+\infty
\end{equation} 
(repeated according to their finite multiplicity) and corresponding eigenfunctions $Y_j(\theta) \in L^2(D)$. Here $\theta=(\theta_1,...,\theta_{N-1})$ is the system of coordinates on $D$ induced by the spherical coordinates in $\R^N$. The eigenfunctions $\{Y_j\}_{j\in\mathbb N}$ can be chosen to form a Hilbert basis for $L^2(D)$,  they satisfy
\beq\label{def:lapleigbe}
-\Delta_{\S^{N-1}}Y_j(\theta)=\lambda_j(D) Y_j(\theta),\quad \theta\in D, 
\eeq
and 
$$
\int_{D} \langle \nabla_{\mathbb S^{N-1}} Y_j(\theta) , \nabla_{\mathbb S^{N-1}} Y_i(\theta)\rangle d\sigma(\theta) = \lambda_j(D)\delta_{ij}
$$
where we denote by $\nabla_{\mathbb S^{N-1}}$ the gradient on $\mathbb S^{N-1}$.

In particular,  $\lambda_0(D)=0$ and the corresponding eigenfunction is constant, while the first non-trivial eigenvalue is $\lambda_1(D)>0$.

We will understand the standard spherical coordinates $\theta$  as polar coordinates in $\mathbb S^{N-1}$: $\theta=(r,\omega)$, where $r\in(0,\pi)$ and $\omega=(\omega_1,...,\omega_{N-2})\in\mathbb S^{N-2}$, where $r$ is the geodesic distance from the north pole.  For $p\in\mathbb S^{N-1}$ we write $p=(r,\omega)$ if $p$ has polar coordinates $(r,\omega)$. 

\subsection{A family of diffeomorphisms on the sphere}
Let $\alpha\in(0,2)$ and consider the map
$$
\phi_{\alpha}:\mathbb S^{N-1}_+\to C_{\alpha}
$$
defined by
$$
\phi_{\alpha}(r,\omega)=(\alpha r,\omega),
$$
where $C_{\alpha}$ is a (spherical) disk centered at the north pole and radius $\frac{\alpha\pi}{2}$. The map $\phi_{\alpha}$ is a diffeomorphism for any $\alpha\in(0,2)$ and it is the identity for $\alpha=1$.

Let now
$$
D_{\alpha}:=\phi_{\alpha}(D)
$$
We clearly have $D_1=D$.

\begin{proposition}\label{limit_L2}

For any smooth $D\subset\S^{N-1}_+$ we have
$$
\lim_{\alpha\to 0}\lambda_j(D_{\alpha})=+\infty \qquad \forall j\ge 1.
$$
\end{proposition}
The proof will be given in Section \ref{ss3.3}. The following corollary immediately follows.
\begin{corollary}\label{cor_L2}
Let $\alpha\in(0,2)$ be such that $D_{\alpha}\subset\S_+^{N-1}$, and $  \lambda_1(D_{\alpha})<N-1$. Then there exists $\bar\alpha\in(0,\alpha)$ such that
$$
\lambda_1(D_{\bar\alpha})=N-1.
$$
Possibly redefining the parameter $\alpha$, we may assume $\bar\alpha=1$.
\end{corollary}

In order to understand the behavior of $\lambda_1(D_{\alpha})$ as $\alpha$ varies, we interpret the eigenvalue $\lambda_1(D_{\alpha})$ as the first Neumann eigenvalue of a suitable operator on the fixed domain $D$. In fact, $\lambda_1(D_{\alpha})$ is the first eigenvalue of the Laplacian on $D$ with respect to the metric $g_{\alpha}$, the pull-back of the round metric on $D_{\alpha}$ through $\phi_{\alpha}$.

More precisely, if $g$ denotes the standard round metric on $\mathbb S^{N-1}$, consider on $D$ the metric $g_{\alpha}:=\phi_{\alpha}^*g$:
$$
\langle X,Y\rangle_{g_{\alpha}}=\langle d\phi_{\alpha}(X),d\phi_{\alpha}(Y)\rangle_g
$$
for all tangent vectors $X,Y$. Here by $\langle\cdot,\cdot\rangle_h$ we denote the scalar product induced on tangent spaces by a metric $h$. In what follows we will highlight the dependence on the metric only when necessary, and we will often omit it when the metric is the  standard Euclidean metric in $\mathbb R^N$ or the standard round metric on $\mathbb S^{N-1}$. 
Then (by construction), the manifold $(D,g_{\alpha})$ is isometric to the manifold $(D_{\alpha},g)$, i.e., the domain $D_{\alpha}$ with the round metric $g$. Therefore
$$
\lambda_1(D_{\alpha})=\lambda_1(D,g_{\alpha})
$$
where the notation $\lambda_1(D,h)$ simply stands for the first Neumann eigenvalue of the Laplacian $\Delta_h$ on $D$ for a metric $h$.
The Laplacian commutes with isometries, that is
$$
\Delta_{g_{\alpha}}(\phi_{\alpha}^*f)=\phi_{\alpha}^*(\Delta_gf)(=\phi_{\alpha}^*(\Delta_{\mathbb S^{N-1}}f))
$$
and this is a concrete way of looking at the new operator on the fixed domain $D$.

If $D\subset\mathbb S^{N-1}_+$, then $r<\pi/2$ for any point of $D$ ($r$ is the polar coordinate, i.e., the distance from the north pole). The metric $g_{\alpha}$ on $D$, for $\alpha\in(0,2)$, is almost isometric to an homothety of the spherical metric on $D$, in fact 
$$
\min\{1,\frac{\sin(\alpha\frac{\pi
}{2} )}{\alpha}\}<\frac{\sin(\alpha r)}{\alpha\sin(r)}<\max\{1,\frac{\sin(\alpha\frac{\pi
}{2} )}{\alpha}\}
$$
and then we write
\beq\label{quasiisom}
g_{\alpha}=\alpha^2dr^2+\sin^2(\alpha r)g_{\mathbb S^{N-2}}\approx\alpha^2(dr^2+\sin^2( r)g_{\mathbb S^{N-2}})=\alpha^2 g.
\eeq

\subsection{The induced diffeomorphisms on the cone}

We pass now to the cone  $\Sigma_{D_{\alpha}}$, which we want to recast to the fixed cone $\Sigma_{D}$. Here, we are considering spherical coordinates in $\mathbb R^N$; in particular, we denote by $\rho$ the distance from the origin in $\mathbb R^N$, that is, the radial coordinate. A point $x\in\Sigma_D$ will have then coordinates $(\rho,p)\in(0,\infty)\times D_{\alpha}$. The standard round metric on $\mathbb S^{N-1}$ is denoted by $g$. The Euclidean metric $G$ on $\mathbb R^N$ (i.e., the metric on $\Sigma_{D_{\alpha}}$) in the coordinates $(\rho,p)$ is
$$
G=d\rho^2+\rho^2 g
$$
Now we consider the map $\Phi_{\alpha}:\Sigma_{D}\to \Sigma_{D_{\alpha}}$  defined by
\beq\label{Phidef}
\Phi_{\alpha}(\rho,p)=(\rho,\phi_{\alpha}(p))
\eeq
and again we consider on $\Sigma_{D}$ the metric $G_{\alpha}$, i.e., the pull-back of the Euclidean metric $G$ though $\Phi_{\alpha}$. Precisely, we have
$$
G_{\alpha}:=\Phi_{\alpha}^*G=d\rho^2+\rho^2\phi_{\alpha}^*g=d\rho^2+\rho^2g_{\alpha}.
$$
Now we express the Laplacian on $\Sigma_{D}$ with respect to the metric $G_{\alpha}$:
\begin{equation}\label{laplacian_new_metric}
\Delta_{G_{\alpha}}=\partial^2_{\rho\rho}u+\frac{N-1}{\rho}\partial_\rho u+\frac{1}{\rho^2}\Delta_{g_{\alpha}}u.
\end{equation}
Again, we have that
\begin{equation}\label{commuting}
\Delta_{G_{\alpha}}(\Phi_{\alpha}^*f)=\Phi_{\alpha}^*(\Delta_Gf)(=\Phi_{\alpha}^*(\Delta f)),
\end{equation}
and $\Delta_G=\Delta$ is the usual Laplacian with respect to the Euclidean metric $G$.

\subsection{Derivative of $\lambda_1(D_{\alpha})$ with respect to $\alpha$}

Now we compute and estimate the derivative
$\frac{d}{d\alpha}\lambda_1(D_{\alpha})\big|_{\alpha=1}=\frac{d}{d\alpha}\lambda_1(D,g_\alpha)\big|_{\alpha=1}$ at a simple eigenvalue. Recall that if $\lambda_1(D_1)$ is simple, then $\lambda_1(D_{\alpha})$ is simple for $\alpha\in(1-\epsilon,1+\epsilon)$ for small $\epsilon$.

\begin{proposition}\label{3.3}
Assume that $\lambda_1(D)=\lambda_1(D_1)$ is simple, and let $u$ be an eigenfunction of $\lambda_1(D)$ normalized by $\int_{D}u^2=1$. Then
\begin{multline}\label{derivative}
\frac{d}{d\alpha}\lambda_1(D_{\alpha})|_{\alpha=1}
=\frac{d}{d\alpha}\lambda_1(D,g_{\alpha})|_{\alpha=1}\\
=\int_{D}\left[(-1+(N-2)\frac{r}{\tan(r)})(\partial_ru)^2+(1+(N-4)\frac{r}{\tan(r)})\frac{|\nabla_{\mathbb S^{N-2}}u|^2}{\sin^2(r)}\right]\sin^{N-2}(r)drd\sigma_{\mathbb S^{N-2}}\\
-\lambda_1(D)\int_{D}u^2(1+(N-2)\frac{r}{\tan(r)})\sin^{N-2}(r)drd\sigma_{\mathbb S^{N-2}}.
\end{multline}
\end{proposition}
To clarify the notation of \eqref{derivative}, here we are considering polar coordinates $(r,\omega)\in(0,\pi/2)\times\mathbb S^{N-2}$ in $\mathbb S^{N-1}_+$ centered at the north pole, and $d\sigma_{\mathbb S^{N-2}}=d\sigma_{\mathbb S^{N-2}}(\omega)$ is the volume element of $\mathbb S^{N-2}$.

The proof of the proposition is standard and consists of deriving the eigenvalue of the Laplacian with respect to a $1$-parameter family of varying metrics, in our case, $g_{\alpha}$, on the fixed manifold $D$. We give the proof following the classical one, which can be found in Berger \cite{berger}.

\begin{proof}[Proof of Proposition \ref{3.3}]
Following \cite[Section 3]{berger}, we find that
$$
\frac{d}{d\alpha}\lambda_1(D,g_{\alpha})|_{\alpha=1}=-\int_{D}\langle q(u),h\rangle_g d\sigma_g,
$$
where $h=\frac{d}{d\alpha}{g_{\alpha}}_{|_{\alpha=1}}$ and $q$ is the quadratic form defined as
$$
q(u)=du\otimes du-\frac{1}{2}(|\nabla u|^2-\lambda_1(D)u^2)g.
$$
The notation $\langle\cdot,\cdot\rangle_g$ stands for the inner product on tensors induced by the metric $g$ and $\otimes$ is the tensor product. We recall that here $g=g_1$ is the standard round metric.

Now we have that 
$$
g=dr^2+\sin^2(r)g_{\mathbb S^{N-2}},
$$
$$
h=2dr^2+2r\sin(r)\cos(r)g_{\mathbb S^{N-2}},
$$
and 
$$
du\otimes du=(\partial_ru)^2dr\otimes dr+2\frac{\partial_ru}{\sin(r)}dr\otimes d_{\mathbb S^{N-2}}u+\frac{d_{\mathbb S^{N-2}}u\otimes d_{\mathbb S^{N-2}}u}{\sin^2(r)}.
$$
All these objects are $(0,2)$ tensors, and then we compute
$$
\langle g,h\rangle_g=2+2(N-2)\frac{r}{\tan(r)}
$$
and
$$
\langle du\otimes du,h\rangle_g=2(\partial_ru)^2+2\frac{r}{\tan(r)}\frac{|\nabla_{\mathbb S^{N-2}}u|^2}{\sin^2(r)}.
$$
Here we have used that $|\nabla_{\mathbb S^{N-2}}u|^2=|d_{\mathbb S^{N-2}}u|^2$. Recalling that $|\nabla u|^2=(\partial_ru)^2+\frac{|\nabla_{\mathbb S^{N-2}}u|^2}{\sin^2(r)}$ and that $d\sigma_g=\sin^{N-2}(r)drd\sigma_{\mathbb S^{N-2}}$ we get formula \eqref{derivative}.
\end{proof}

Our aim is to show that the derivative \eqref{derivative} is negative under suitable assumption on $D$. 

First, note that $\frac{r}{\tan(r)}\in(0,1)$ when $r\in(0,\pi/2)$ ($\frac{r}{\tan(r)}\to 0$ as $r\to\pi/2$, $\frac{r}{\tan(r)}\to 1$ as $r\to 0$, and it is decreasing in $r$). 

\medskip

\begin{proposition}
    Let $\lambda_1(D)=\lambda_1(D_1)$ be simple and assume that $D\subset C_t$ where $C_t$ is a geodesic disk of radius $t\in(0,\pi/2)$ such that $t/\tan (t)\ge\frac{N-4}{N-2}$. Then
\beq\label{derneg}
\frac{d
}{d\alpha}\lambda_1(D_\alpha)|_{\alpha=1}<0.
\eeq
\end{proposition}

\begin{proof}
    From \eqref{derivative} and the fact that $\frac{r
    }{\tan r}<1$ for $r\in(0,\pi/2)$ we have
\begin{multline}\label{derivative2}
\frac{d}{d\alpha}\lambda_1(D_{\alpha})|_{\alpha=1}
\\
<\int_{D}\left[(N-3)\left((\partial_ru)^2+\frac{|\nabla_{\mathbb S^{N-2}}u|^2}{\sin^2(r)}\right)-\lambda_1(D)u^2((N-2)b+1)\right]\sin^{N-2}(r)drd\sigma_{\mathbb S^{N-2}}\\
=(N-3)\int_{D}|\nabla u|^2-\left(\frac{(N-2)b+1}{N-3}\right)\lambda_1(D)u^2\sin^{N-2}(r)drd\sigma_{\mathbb S^{N-2}},
\end{multline}
    where $b=b(r)=r/\tan r$.
Then taking $t$ such that $t/\tan t\ge\frac{N-4}{N-2}$ by the monotonicity of the function $t/\tan t$ we have that
\beq\label{contob}
b=\frac{r}{\tan r}\ge \frac{t}{\tan t}\ge \frac{N-4}{N-2}\quad \forall r\in (0,t),
\eeq
and hence
$$
\frac{(N-2)b+1}{N-3}>1.
$$
Since 
$$
\lambda_1(D)=\frac{\int_D |\nabla_{\mathbb S^{N-1}} u|^2}{\int_D u^2},
$$
we obtain \eqref{derneg}.
    
\end{proof}
\begin{remark}
Note that for $N=3,4$ one can see directly from \eqref{derivative2} or from \eqref{contob} that the derivative is negative for any  domain $D_1=D\subset\mathbb S^{N-1}_+$. Instead for $N\geq 5$ we have to consider domains contained in smaller disks.
\end{remark}

\subsection{Proof of Proposition \ref{limit_L2}}\label{ss3.3}
We prove here Proposition \ref{limit_L2}. In fact, one can see that $g_{\alpha}$ is quasi-isometric to a homothety of the round metric $g$ of factor $\alpha^2$, and the eigenvalues can be compared.

\begin{proof}[Proof of Proposition \ref{limit_L2}]

Observe that from \eqref{quasiisom} on $D$ the metric $g_{\alpha}$ and the metric $\alpha^2 g$ are quasi-isometric. This means that there exists a constant $K>0$ such that
\beq\label{compar}
\frac{1}{K}<\frac{g_{\alpha}}{\alpha^2 g}< K
\eeq
For any eigenvalue $\lambda_j$, by the min-max principle we have 
$$
\lambda_j(D,g_\alpha) = \min_{\substack{V\subset H^1(D)\\\dim V = j}} \max_{u \in V, u \neq 0} \frac{\int_{D} |\nabla_{g_{\alpha}} u|_{g_{\alpha}}^2 \, d\sigma_{g_{\alpha}}}{\int_{D} u^2 \, d\sigma_{g_{\alpha}}}.
$$
Here $d\sigma_{g_{\alpha}}$ is the Riemannian volume element of the metric $g_{\alpha}$ and $|\nabla_{g_{\alpha}}u|_{g_{\alpha}}$ is the norm of the gradient of $u$ for the metric $g_{\alpha}$.
From \eqref{compar} we deduce
$$
d\sigma_{\alpha^2g} = \sqrt{\det \alpha^2g} \, d\theta_1 \dots d\theta_{N-1} > \frac{1}{K^{\frac{N-1}{2}}} \sqrt{\det g_\alpha} \, d\theta_1 \dots d\theta_{N-1} = \frac{1}{K^{\frac{N-1}{2}}} d\sigma_{g_\alpha}
$$
and
$$
d\sigma_{\alpha^2g}<K^{\frac{N-1}{2}} d\sigma_{g_\alpha}.
$$
Moreover,
$$
\frac1K|\nabla_{g_\alpha} u|^2_{g_\alpha}<|\nabla_{\alpha^2g} u|^2_{\alpha^2g}<K|\nabla_{g_\alpha} u|^2_{g_\alpha}.
$$
Then 

$$\frac{\int_{D} |\nabla_{g_{\alpha}} u|_{g_{\alpha}}^2 \, d\sigma_{g_{\alpha}}}{\int_{D} u^2 \, d\sigma_{g_{\alpha}}}<K^{1+(N-1)/2+(N-1)/2}\frac{\int_{D} |\nabla_{\alpha^2g} u|_{\alpha^2g}^2 \, d\sigma_{\alpha^2g}}{\int_{D} u^2 \, d\sigma_{\alpha^2g}}
$$

and
$$\frac{1}{K^{1+(N-1)/2+(N-1)/2}}\frac{\int_{D} |\nabla_{\alpha^2g} u|_{\alpha^2g}^2 \, d\sigma_{\alpha^2g}}{\int_{D} u^2 \, d\sigma_{\alpha^2g}}<\frac{\int_{D} |\nabla_{g_{\alpha}} u|_{g_{\alpha}}^2 \, d\sigma_{g_{\alpha}}}{\int_{D} u^2 \, d\sigma_{g_{\alpha}}}.
$$
Then, from the Min-Max Theorem, we conclude that
$$
\frac{1}{K^N}\lambda_j(D,\alpha^2g)<\lambda_j(D,g_\alpha)<K^N\lambda_j(D,\alpha^2g).
$$

Now, we have that
$$
\lambda_j(D,\alpha^2 g)=\frac{1}{\alpha^2}\lambda_j(D,g)=\frac{1}{\alpha^2}\lambda_j(D)
$$
and this yields the result.

\end{proof}

\subsection{An example of simple eigenvalue $N-1$}

\begin{figure}
\includegraphics[width=0.9\textwidth]{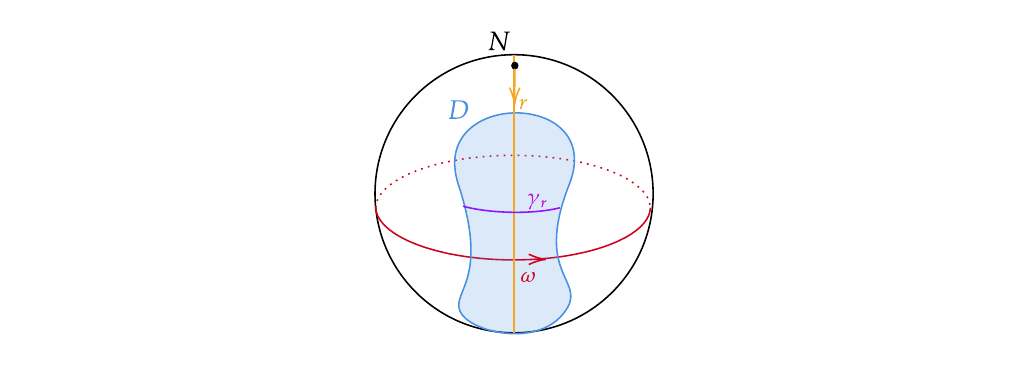}
\caption{}
\label{fig1}
\end{figure}

We provide here an example of a domain $D$ contained in an hemisphere in $\R^3$ for which the first eigenvalue $\lambda_1(D)$ is $2$ and is simple. We just start from a dumbbell domain on $\mathbb S^2_+$. Assume it is symmetric with respect to the equator and to the (arc of) great circle $\omega=0$. Here we are using polar coordinates $( r ,\omega)$, with $ r $ the distance from the north pole. See Figure \ref{fig1}. Roughly speaking, $D$ consists of two disjoint `fat' parts (e.g., two geodesic disks) connected by a thin channel. In our example, the fat parts are symmetric with respect to the equator (see Figure \ref{fig1}). The domain $D$ is contained in the strip $ r \in( r _0,\pi- r _0)$, $ r _0\in(0,\pi/2)$. For any $ r '\in( r _0,\pi- r _0)$, we call $\gamma_{ r '}:=D\cap\{ r = r '\}$ (see Figure \ref{fig1}). Assume that $\max_{ r \in( r _0,\pi- r _0)}|\gamma_{ r }|<\pi/2$.

From Corollary \ref{cor_L2} and from the fact that the first eigenvalue of a fixed dumbbell can be arbitrarily small when the width of the channel is small (see e.g., \cite{IPW}), we deduce that we can realize such a domain $D=D_1$ with $\lambda_1(D_1)=2$. We show that $\lambda_1(D_1)=2$ it is simple. It will follow that $\lambda_1(D_{\alpha})$  remains simple for all $\alpha\in(1-\varepsilon,1+\varepsilon)$ for some $\varepsilon>0$.

 Since we are dealing with $\lambda_1(D)$, its multiplicity is at most $2$ by Courant's Theorem.  Since any corresponding eigenfunction has exactly two nodal domains, by symmetry  we know that the eigenspace is spanned by eigenfunctions that are either odd with respect to the equator and even with respect to the arc $\omega=0$, or the other way around. Note that an eigenfunction $u$ corresponding to $\lambda_1(D)$ cannot be odd (respectively even) with respect to both the equator and the great circle $\omega=0$. Assume that $u$ is odd with respect to the arc of great circle $\omega=0$ (the yellow arc in Figure \ref{fig1}). Then, for any $ r \in( r _0,\pi- r _0)$, $\int_{\gamma_{ r }}u=0$. Therefore, $\int_{\gamma_{ r }}(\partial_{\omega}u)^2\geq\frac{\pi^2}{|\gamma_{ r }|^2}\int_{\gamma_{ r }}u^2$. Here $\frac{\pi^2}{|\gamma_{ r }|^2}$ is the first (non trivial) Neumann eigenvalue on $\gamma_{ r }$.

By assumption we have $|\gamma_{ r }|<\pi/2$. This implies that $\int_{\gamma_{ r }}(\partial_{\omega}u)^2\geq 4\int_{\gamma_{ r }}u^2$, which means $\int_{\gamma_{ r }}|\nabla u|^2\geq 4\int_{\gamma_{ r }}u^2$. Integrating this last inequality in $ r \in( r _0,\pi- r _0)$ we get that $\int_D|\nabla u|^2\geq 4\int_Du^2$. Hence we have a contradiction since $\lambda_1(D)=2$. Therefore the first eigenvalue is simple and a corresponding eigenfunction is odd with respect to the equator.

The construction is valid in any space dimension, we just needed to consider a rotationally symmetric dumbbell in $\mathbb S^{N-1}$ foliated by $N-2$ dimensional spheres of small enough maximal volume.

\section{Bifurcation results}\label{sec3}
\subsection{Spectral analysis}\label{spectan}
Let $D\subset\S^{N-1}_+$ be a smooth domain and assume that the corresponding cone $\Sigma_D$ is a Lipschitz domain. We consider the map $\phi_\alpha, \alpha \in (0,2)$ defined in Section \ref{ss3.1} and the cone $\Sigma_{D_\alpha}$.
Consider the bubble $U\in\mathcal{
D
}^{1,2}(\Sigma_{D_\alpha})$:
\beq\label{bubbledef}
U(x):=c_0(1+|x|^2)^{-\frac{N-2}{2}}\qquad x\in \Sigma_{D_\alpha},
\eeq
where $c_0$ is a constant such that 
$$
-U^{2-2^*}\Delta U=U\qquad \text{in} \quad \Sigma_{D_\alpha}.
$$
The rescaled bubbles will be denoted by $$
U_{s}(x):=s U(s^{\frac{2}{N-2}}x).
$$
\begin{lemma}\label{lemma_spettro}
The operator $L_{U}=U^{2-2^*}\Delta$ on $L^2(U^{2^*-2})$ has discrete spectrum.
\end{lemma}
\begin{proof}
    The proof can be found in \cite{CPP2}.
\end{proof}
To describe the spectrum of the operator $L_{U}$ we need to know the eigenvalues $\lambda_j(D_\alpha), j\in \N$, of problem \eqref{pb_lambda_def1} on $D_\alpha$.
In particular, we recall that the  $\lambda_0(D_\alpha)=0$ and the corresponding eigenfunction is constant, while the first positive eigenvalue is $\lambda_1(D_\alpha)$.

Now we are ready to address the following eigenvalue problem that is studied in \cite{CPP2},
\begin{equation}\label{pb_mu_def2}
\begin{cases}
\Delta v+\mu_{\alpha} U^{2^*-2}v=0 &\text{in } {\Sigma_{D_\alpha}}\\
\partial_\nu v=0 &\text{on }  \partial{\Sigma_{D_\alpha}}\setminus\{0\} \\
v\in \mathcal{D}^{1,2}(\Sigma_{D_\alpha}), &
\end{cases}
\end{equation}
whose  eigenvalue will be denoted by $\mu_{i,\alpha}, i\in\N_+$. The eigenvalue equation in \eqref{pb_mu_def2} can be written by using spherical coordinates, and we can look for solutions by using separation of variables: let $(\rho,\theta)\in(0,+\infty)\times\mathbb S^{N-1}$ the standard spherical coordinates in $\mathbb R^N$ and let $v=R(\rho)Y(\theta)$ be a solution of \eqref{pb_mu_def2} for an eigenvalue $\mu_\alpha$, then it holds that $Y$ is an eigenfunction of \eqref{pb_lambda_def1} on $D_\alpha$ for some $j$, while $R$ satisfies (see \cite{CPP2}, Proposition 2.6)
\begin{equation} \label{SL_eq}
\begin{cases}
(\rho^{N-1}R')'+\rho^{N-1}(-\lambda_j(D_\alpha) \rho^{-2}+\mu_\alpha k_0^{2^*-2}(1+\rho^2)^{-2})R=0 & \text{in }\R_+ \\
\lim_{\rho\to 0}\rho^{N-1}R'(\rho)=0 & \\ 
R(\rho)=o(\rho^{2-N+\epsilon})  \quad \forall\epsilon> 0\text{ as }\rho\to + \infty \,. &

\end{cases}
\end{equation}
We denote by $\mu_{k,\lambda_j(D_\alpha)}, k\in\N_+, j\in\N$, the sequence of eigenvalues of \eqref{SL_eq} obtained in correspondence of $\lambda_j(D_\alpha)$ and denote by $R^j_k$ the corresponding eigenfunctions.

For $\lambda_0(D_\alpha)=0$ the eigenvalues are given by
\beq\label{lambdazerok}\mu_{k,\lambda_0}=\bigg((k-1)(k+N-2)\frac{4}{N(N-2)}+1\bigg)\,.
\eeq
The eigenvalues associated with $\lambda_1(D_\alpha)$ are
\begin{equation}\label{betaeigenvalue}
\mu_{1,\lambda_1}= \frac{1}{N(N-2)}\sqrt{(N-2)^2 + 4\lambda_1(D_\alpha)}( 2 +\sqrt{(N-2)^2 + 4\lambda_1(D_\alpha)}) .  
\end{equation}
In details, it holds the following result (see \cite{CPP2}).

\begin{proposition} \label{lemma_autov}

Let $D_\alpha\subset\S_+^{N-1}$ and $m\in\N_+$ be such that
\beq\label{lemma_autov_cond}
(m-1)(m+N-3)<\lambda_1(D_\alpha)\le m(m+N-2).
\eeq
Then each eigenvalue $\mu_{i,\alpha}$ of \eqref{pb_mu_def2}, with $i=1,...,m$,
is simple and is given by
$$\mu_{i,\alpha}=\bigg((i-1)(i+N-2)\frac{4}{N(N-2)}+1\bigg),
$$
i.e. are the eigenvalues corresponding to $\lambda_0$ in \eqref{SL_eq} while
$$
\mu_{m+1,\alpha}=\frac{1}{N(N-2)}\sqrt{(N-2)^2 + 4\lambda_1(D_\alpha)}\bigg( 2 +\sqrt{(N-2)^2 + 4\lambda_1(D_\alpha)}\bigg).
$$
Moreover, $\psi_i=R_i(r)Y_0(\theta)$ is an eigenfunction of \eqref{pb_mu_def2} corresponding to $\mu_{i,\alpha}$, for any $i=1,...,m,$ where
$$
R_i=\frac{P_{i}(|x|^2)}{(1+|x|^2)^{\frac{N-2}{2}}}
$$
and $P_i(|x|^2)$ are given as in \cite[Appendix]{CPP2}. Instead an eigenfunction corresponding to $\mu_{m+1,\alpha}$ is given by  $\psi_{m+1}=R_{m+1}Y_1(\theta)$ where
$$R_{m+1}=(|x|^2+1)^{1-\frac{N}{2}-\beta(\lambda_1(D_\alpha))}|x|^{\beta(\lambda_1(D_\alpha))}$$
with
$2\beta(\lambda_1(D_\alpha))=-(N-2) + \sqrt{(N-2)^2 + 4\lambda_1(D_\alpha)} $ and $Y_1$ is an eigenfunction of \eqref{pb_lambda_def1} corresponding to $\lambda_1(D_\alpha)$.
\end{proposition}
\begin{proof}
    The proof follows by \cite[Proposition 2.7]{CPP2} and the fact that we are taking $U$ such that $-\Delta U=U^{2^*-1}$, so that the constant $S^{2^*}_U$ in \cite{CPP2} is replaced by $1$.
\end{proof}

\begin{corollary}\label{cor_autoval}
    Let $D_{\bar{\alpha}}\subset\S_+^{N-1}$ be such that $\lambda_1(D_{\bar\alpha})=N-1$, then
    $$
    \mu_{1,{\bar\alpha}}=1
    $$
    and the first eigenfunction is given by $U_{\bar\alpha}
    $,
$$\mu_{2,{\bar\alpha}}= \frac{1}{N(N-2)}\sqrt{(N-2)^2 + 4\lambda_1(D_{\bar\alpha})}\bigg( 2 +\sqrt{(N-2)^2 + 4\lambda_1(D_{\bar\alpha})}\bigg)=
2^*-1
$$
and the second eigenfunctions are given by $\partial_s U_{s,{\bar\alpha}}$ and $\psi_2=(|x|^2+1)^{-\frac{N}{2}}|x|Y_1(\frac{x}{|x|})$. 
\end{corollary}

\begin{proof}
    See \cite{CPP2} Corollary 2.8.
\end{proof}
From now on, as in Corollary \ref{cor_L2} we assume $\bar\alpha=1$, by redefining the parameter $\alpha$, so that $U_{\bar \alpha}$ will be denoted simply by $U$.
\begin{remark}
We say that the bubble $U$ is nondegenerate if the second variation $Q_U$ of the Sobolev functional vanishes only on the 2-dimensional subspace spanned by $\{U,\partial_sU_s|s=1\}$ (see Definition 2.4 of \cite{CPP2}). It holds that the bubbles are nondegenerate if and only if $\lambda_1(D_1)>N-1$ (see Remark 2.5 of \cite{CPP2}and \cite{CPP}). Obviously from Corollary \ref{cor_autoval} for $\lambda_1(D_1)=N-1$ we are in the case where the bubbles are degenerate.
\end{remark}

\subsection{Main results}
Let $D_\alpha$ be a domain as in Section \ref{sec2}.
We consider the following problem
\begin{equation}\label{problem}
\begin{cases}
    -\Delta v=v^{2^*-1} &\text{ in }\Sigma_{D_\alpha}\\
    v\in\mathcal{D}^{1,2}(\Sigma_{D_\alpha}).
\end{cases}
\end{equation}
For the reader convenience we recall the notation for the diffeomorphism considered in Section \ref{sec2}.
 Let $\alpha\in(0,\alpha^*)$, where $\alpha^*$ is the maximum value of $\alpha$ such that $\Phi_\alpha(D_1)\subset \S^{N-1}_+$, and consider the map
$$
\phi_{\alpha}:\mathbb S^{N-1}_+\to C_{\alpha}
$$
defined by
$$
\phi_{\alpha}(r,\omega)=(\alpha r,\omega),
$$
where $C_{\alpha}$ is a (spherical) disk centered at the north pole and radius $\alpha\pi/2$. Then we define
$$
\Phi_{\alpha}(\rho,p)=(\rho,\phi_{\alpha}(p)),
$$
where $(\rho,p)\in(0,+\infty)\times D$, and $p=(r,\omega)\in D$ is a point in $D$. The map $\phi_{\alpha}$ induces
$$
\Phi_{\alpha}^*:\mathcal{D}^{1,2}(\Sigma _{D_\alpha})\to\mathcal{D}^{1,2}(\Sigma _{D_1})
$$
defined by $\Phi_{\alpha}^*u=u\circ \Phi_{\alpha}$. We define an operator $\mathcal L_{\alpha}$ acting on $\mathcal D^{1,2}(\Sigma_{D_1})$ by
$$
\mathcal L_{\alpha}(\Phi_{\alpha}^*v)=\Phi_{\alpha}^*(-\Delta v),
$$
for $v\in\mathcal D^{1,2}(\Sigma_{D_{\alpha}})$ (as usual, the operator is defined in the weak sense). This operator can be described more explicitly: it is $-\Delta_{G_{\alpha}}$ defined in \eqref{laplacian_new_metric}, however we won't need this explicit expression in the following.

In particular, if $v$ is a solution of \eqref{problem}, setting $w=\Phi^*_{\alpha}v$,  we have that $w$ is a weak solution of
\begin{equation}\label{PBTRA}
    \begin{cases}
        \mathcal{L}_\alpha(w)=w^{2^*-1}&\text{ in }\Sigma_{{D_1}}\\
        w \in\mathcal{D}^{1,2}(\Sigma_{D_1}).\\
    \end{cases}
\end{equation}

In fact, since $v$ is a weak solution of \eqref{problem} for every $\varphi\in\mathcal{D}^{1,2}(\Sigma_{D_\alpha})$ it holds
$$
\int_{\Sigma_{D_\alpha}} \langle\nabla v , \nabla\varphi \rangle=\int_{\Sigma_{D_\alpha}}v^{2^*-1}\varphi.
$$
Considering the change of variables in the following way $y=\Phi_{\alpha}(x)$, $x\in D_1$, and setting
$$
w=\Phi_{\alpha}^*v\,,\ \ \ \psi=\Phi_{\alpha}^*\varphi\,,
$$
we get:
$$
\int_{\Sigma_{D_{\alpha}}}\langle\nabla v,\nabla\varphi\rangle=\int_{\Sigma_{D_1}}\langle\nabla_{G_{\alpha}}w,\nabla_{G_{\alpha}}\psi\rangle_{G_{\alpha}} d\sigma_{G_{\alpha}}
$$
and
$$
\int_{\Sigma_{D_{\alpha}}}v^{2^*-1}\varphi=\int_{\Sigma_{D_1}}w^{2^*-1}\psi d\sigma_{G_{\alpha}},
$$
where $G_{\alpha}=\Phi_{\alpha}^*(G)$ is the pull-back of the Euclidean metric $G$ on $\Sigma_{D_{\alpha}}$ through $\Phi_{\alpha}$ and $d\sigma_{G_{\alpha}}$ is the corresponding volume element.

Hence $w$ solves
$$
\int_{\Sigma_{D_1}}\langle\nabla_{G_{\alpha}}w,\nabla_{G_{\alpha}}\psi\rangle_{G_{\alpha}} dv_{G_{\alpha}}=\int_{\Sigma_{D_1}}w^{2^*-1}\psi d\sigma_{G_{\alpha}}
$$
for all $\psi\in\mathcal D^{1,2}(\Sigma_{D_1})$, and this is
 the weak formulation of problem \eqref{PBTRA}.
 
 \begin{remark}
To be precise, the function space for the weak formulation of \eqref{PBTRA} is the space of functions in $L^{2^*}(\Sigma_{D_1},d\sigma_{G_{\alpha}})$ such that $|\nabla_{G_{\alpha}}u|_{G_{\alpha}}\in L^2(\Sigma_{D_1},d\sigma_{G_{\alpha}})$. However, since for all $\alpha$'s we are considering  $G_{\alpha}$ and $G$ are quasi-isometric, this space is equivalent to $\mathcal D^{1,2}(\Sigma_{D_1})$.
\end{remark}

We define the space
$$
X=\big\{u\in \mathcal{D}^{1,2}(\Sigma_{D_1})\big| \sup_{x\in \Sigma_D}\frac{|u(x)|}{U(x)}<\infty\big\},
$$

Then $X$ is a Banach space equipped with the norm $\|u\|_X:=\max\big\{\|u\|_{1,2};\sup_{x\in\Sigma_D}\frac{|u(x)|}{U(x)}\big\}$. Observe that $\mathcal{L_\alpha}:X\to L^2(U^{2^*-2})$.

To rule out the degeneracy due to the invariance under dilation of problem \eqref{problem}, we will solve the linearized equation in the subspace of functions that are invariant under the Kelvin transform $k(v)$ which is defined by
$$
k(u)(x):=\frac{1}{|x|^{N-2}}u\bigg(\frac{x}{|x|^2}\bigg)
$$
and we denote by $K$ the subspace of functions in $X$ which are invariant by the Kelvin transform, i.e.
$$
K:=\{v\in X|k(v)=v \}.
$$
Observe that
$k(U)=U$ and $k(\partial_sU_{s}|_{s=1})=-\partial_sU_{s}|_{s=1}$.

\begin{remark}\label{rmkpalla}
Let us consider a ball $B:=\{u\in X|\|u-U\|_X<\delta\}$ and the set $K\cap B$. If $\delta$ is sufficiently small then any $u\in K\cap B$ is positive on $\bar \Sigma_D$.
Indeed if this was not true there would exist a sequence of functions $u_{n}\in X\cap K$ and points $x_n\in \bar\Sigma_D$ such that
$$
u_n\to U \quad \text{ in }X \quad\text{  and  }\quad u_{n}(x_n)\le 0.
$$
We have two cases:

Case 1) $x_n\to \bar x\in \bar\Sigma_D$.  Then it holds
\beq\label{conto}
|u_{n}(x_n)-U(\bar x)|\le|u_{n}(x_n)-U(x_n)|+|U(x_n)-U(\bar x)|\le \|u_{n}-U\|_\infty+|U(x_n)-U(\bar x)|\to 0 ,
\eeq

as $n\to\infty$, since $u_{n}\to U$ in $X$ and $U$ is continuous. 
From \eqref{conto} and the fact that $U>0$ on $\bar\Sigma_D$ we get a contradiction. 

Case 2) $x_n\to \infty$. In this case, since $u_{n}$ and $U$ are Kelvin invariant we apply the same argument as before to the Kelvin transform of $u_{n}$. 

\end{remark}

\begin{definition}\label{Tdef}
    Let us consider the ball $B:=\{u\in X|\|u-U\|_X<\delta\}$ where $\delta$ is defined as in Remark \ref{rmkpalla}, and define the operator
    $$
    F:(0,\alpha^*)\times (K\cap B)\to K
    $$
\beq\label{defoperator}
    F(\alpha, u):= u-( \mathcal{L}_\alpha
    )^{-1}(u^{2^*-1}).
    \eeq
    \end{definition}
   
We want to prove the bifurcation result using the following well known Theorem \cite[Theorem 1.7]{CR} applied to the operator $F$.
\begin{theorem}[Crandall-Rabinowitz Bifurcation Theorem] \label{CRT}
Let $X$ and $Y$ be Banach spaces, and let $\Omega\subset X$ and $I \subset \mathbb{R}$ be open domains, where we assume $0 \in \Omega$. Denote the elements of $\Omega$ by $v$ and the elements of $I$ by $t$. Let $F: I \times \Omega \to Y$ be a $C^2$ operator such that

\begin{itemize}
    \item[(i)] $F(t, 0) = 0$ for all $t \in I$,
    \item[(ii)] $\text{Ker} D_v F(t_*, 0) = \mathbb{R}w$ for some $t_* \in I$ and some $w \in X \setminus \{0\}$;
    \item[(iii)] $\text{codim} \, \text{Im} D_v F(t_*, 0) = 1$;
    \item[(iv)] $D_t D_v F(t_*, 0)(w) \notin \text{Im} D_v F(t_*, 0)$.
\end{itemize}

Then there exists a nontrivial $C^1$ curve
$$\quad (-\epsilon, \epsilon) \ni s \mapsto (t(s), v(s)) \in I \times X,
$$
for some $\epsilon > 0$, such that:
\begin{itemize}
    \item[(1)] $t(0) = t_*, \, t'(0) = 0, \, v(0) = 0, \, v'(0) = w$;
    \item[(2)] $F(t(s), v(s)) = 0$ for all $s \in (-\epsilon, +\epsilon)$.
\end{itemize}

Moreover, there exists a neighborhood $N$ of $(t_*, 0)$ in $X \times Y$ such that all solutions of the equation $F(t, v) = 0$ in $N$ belong to the trivial solution line $\{(t, 0)\}$ or to the curve. The intersection $(t_*, 0)$ is called a bifurcation point.

\end{theorem}
\begin{lemma}
 Let $F$ be as defined in Definition \ref{defoperator}, the zeros of the operator $F$ satisfy
$$
\begin{cases}
     \mathcal{L}_\alpha u=u^{2^*-1} &\Sigma_{D_1}\\

    u\in K\cap B,
\end{cases}$$
in particular $F(\alpha, U)=0$ for all $\alpha\in(0,\alpha^*)$.
\end{lemma}
\begin{proof}
The zeros of the operator $F$ satisfy
$$
u-( \mathcal{L}_\alpha
)^{-1}(u^{2^*-1})=0
$$
namely
$$
\begin{cases}
     \mathcal{L}_\alpha u=u^{2^*-1} &\Sigma_{D_1}\\

    u\in K\cap B.
\end{cases}$$

\end{proof}

\begin{lemma}\label{contder}
    The operator $F$ is continuous from $(0,\alpha^*)\times (K\cap B)$ into $ K$ and its derivatives $\partial_u F, $ $\partial_{\alpha}(\partial_uF)$ exist and are continuous.
\end{lemma}
\begin{proof}
Clearly  $F:(0,\alpha^*) \times (K\cap B)\mapsto K$, because it is a composition of operators that map $(0,\alpha^*) \times (K\cap B)\mapsto K$. Thanks to the critical exponent and the invariance of the Laplacian and of the diffeomorphisms under Kelvin transform, we can conclude that 
 $F:(0,\alpha^*)\times (K\cap B)\mapsto K$.
 To prove the continuity of $F$ in $K\cap B$, let $\alpha_n\to\alpha$ in $\R$ and $u_n\to u$ in $X$ as $n\to \infty$.
Since $u_n\to u$ in $\mathcal{D}^{1,2}(\Sigma_{D_1})$,  the convergence also holds in $L^{2^*}(\Sigma_{D_1})$. 
We write
\[
F(\alpha_n,u_n)-F(\alpha,u)
=(u_n-u)-
\Big[(\mathcal L_{\alpha_n})^{-1}(u_n^{2^*-1})-(\mathcal L_{\alpha})^{-1}(u^{2^*-1})\Big].
\]
The first term $u_n-u\to 0$ in $X$ by assumption.
For the second term we split:
\[
(\mathcal L_{\alpha_n})^{-1}(u_n^{2^*-1})-(\mathcal L_{\alpha})^{-1}(u^{2^*-1})
= (I)+(II),
\]
where
\[
(I)=\big((\mathcal L_{\alpha_n})^{-1}-(\mathcal L_{\alpha})^{-1}\big)(u^{2^*-1}), \quad
(II)=(\mathcal L_{\alpha_n})^{-1}(u_n^{2^*-1}-u^{2^*-1}).\qquad
\]

\smallskip
Convergence of (I).
Since the operator $\mathcal L_\alpha$ is continuous in operator norm we have easily that
\[
\| (I)\|_{\mathcal D^{1,2}} \to 0.
\]
Hence $(I)\to 0$ in $X$.

\smallskip
Convergence of (II).
Since $u_n\to u$ in $L^{2^*}$ and $|u_n|\le C U$ uniformly, we obtain
\[
u_n^{2^*-1}\to u^{2^*-1}\quad\text{in }L^{2^*}.
\]
Again by continuity of $(\mathcal L_{\alpha_n})^{-1}$ as an operator
$L^{2^*}\to\mathcal D^{1,2}$, we deduce
\[
\| (II)\|_{\mathcal D^{1,2}}\to 0,
\]
and thus $(II)\to 0$ also in $X$.

\smallskip
Combining the two contributions we conclude
\[
F(\alpha_n,u_n)\to F(\alpha,u)\quad\text{in }X.
\]

The existence of $\partial_u F$, and $\partial_\alpha(\partial_uF)$ (for the topology of X) and their continuity follows in a similar way and we omit it.
\end{proof}
\begin{lemma}\label{ker}
     Let $F$ be as defined in \eqref{defoperator}, and assume that $\lambda_1({D_1})$ is simple. Then $ker(\partial_u F(1,U))$ is one dimensional and it is given by

     $$
     ker(\partial_u F(1,U))=span\{\psi_2\},
     $$
     where $\psi_2$ is the eigenfunction defined in Corollary \eqref{cor_autoval}.
\end{lemma}
\begin{proof}
     Let us consider the Fréchet derivative of $F$ at $(1,U)$. We have that
  \beq\label{deroper}
    \partial_u F(\alpha,U)[w]=w-(2^*-1)( \mathcal{L}_\alpha)^{-1}[U^{2^*-2}w], \quad w\in K
\eeq
    so that $\partial_u F(\alpha,U)(w)=0$ if and only if $w\in K$ is a solution to
    $$
    \mathcal{L}_\alpha w-(2^*-1)U^{2^*-2}w=0.
    $$
    From Proposition \ref{lemma_autov} and Corollary \ref{cor_autoval}, it follows that the only solutions to the equation for $\alpha=1$ are given by
    $\partial_sU_s$ and $\psi_2:=R_2(|x|)Y_1(\frac{x}{|x|})$, where $R_2(x)=(1+|x|^2)^{-\frac{N}{2}}|x|$ and $Y_1$ is the only eigenfunction of $\lambda_1({D_1})=N-1$ (we are supposing that $\lambda_1(D_1)$ is simple). Now $k(\partial_sU_s)=-\partial_sU_s$, so it is not invariant under Kelvin transform, while
    $$
k(\psi_2)(x)=|x|^{2-N} (1+|x|^{-2})^{-\frac{N}{2}}|x|^{-1}Y_1(\frac{x}{|x|})=\psi_2
    $$
    so $\psi_2\in K$. Then $w\in  ker(\partial_u F(1,U))$ if and only if $w\in span\{\psi_2\}$.
\end{proof}

\begin{lemma}\label{range}
    Let $F$ be as defined in \eqref{Tdef}, the range $Ran(\partial_uF(1,U))\subset K$ has codimension one. It is the set of the functions orthogonal to $\psi_2$.
\end{lemma}
\begin{proof}
   Since 
    $$
\partial_u F(\alpha,U)[w]=w-(2^*-1)( \mathcal{L}_\alpha)^{-1}U^{2^*-2}w, \qquad w\in K
$$
by the compactness of the embedding $\mathcal{D}^{1,2}(\Sigma_D)\hookrightarrow L^2(U^{2^*-2})$ (see \cite{CPP2}) we have that $\partial_u F(1,U)$ is a compact perturbation of the identity. Thus, the thesis follows from the Fredholm Alternative.
\end{proof}

\begin{lemma}\label{transvers}
     Let $F$ be defined by \eqref{defoperator} and $D_1\subset C_t$ where $C_t$ is a geodesic disk of radius $t\in(0,\pi/2)$ such that $t=\frac{N-4}{N-2}\tan(t)$ ($t=\pi/2$ if $N=3,4$). Then we have

     $$
     \partial_{\alpha }(\partial_uF(\alpha,U)[\psi_2])\big|_{\alpha=1}\not\in Ran(\partial_u F(1,U)).
     $$
\end{lemma}
\begin{proof}
    From \eqref{deroper} we have
    $$
    \partial_u F(\alpha,U)[\psi_2]=\psi_2-(2^*-1)( \mathcal{L}_\alpha)^{-1}[U^{2^*-2}\psi_2], 
    $$
    then differentiating with respect to $\alpha$ yields
    $$
    \partial_{\alpha }(\partial_uF(\alpha,U)[\psi_2])=-(2^*-1) \frac{d}{d\alpha}( \mathcal{L}_\alpha)^{-1}[U^{2^*-2}\psi_2]. 
    $$

By Lemma \ref{range} we have to show that $(- \frac{d}{d\alpha}( \mathcal{L}_\alpha)^{-1}[U^{2^*-2}\psi_2]\big|_{\alpha=1},\psi_2)\ne 0$.  Thus, we have to show that
$$
\int_{\Sigma_{D_1}}-\langle\nabla( \frac{d}{d\alpha}( \mathcal{L}_\alpha)^{-1}[U^{2^*-2}\psi_2]\big|_{\alpha=1}),\nabla\psi_2\rangle\ne 0
$$
Hence at $\alpha=1$, from the eigenvalue problem \eqref{pb_mu_def2} and Corollary \ref{cor_autoval} we have that on $\Sigma_{D_1}$
$$
( \mathcal{L}_\alpha)^{-1}[U^{2^*-2}\psi_2]=\frac{1}{ \mu_{2,\alpha}}\psi_2
$$
and hence
$$
 \frac{d}{d\alpha}( \mathcal{L}_\alpha)^{-1}[U^{2^*-2}\psi_2]\big|_{\alpha=1}=-\frac{1}{(\mu_{2,\alpha})^2} \frac{d}{d\alpha} \mu_{2,\alpha}|_{\alpha=1}\psi_2.
$$
Then we need to show that
$$
 \frac{d}{d\alpha}\mu_{2,\alpha}|_{\alpha=1}\int_{\Sigma_{D_1}}|\nabla\psi_2|^2\ne 0.
$$
 From Proposition \ref{lemma_autov} we  know that
    $$\mu_{2,{\alpha}}= \frac{1}{N(N-2)}\sqrt{(N-2)^2 + 4\lambda_1(D_{\alpha})}\bigg( 2 +\sqrt{(N-2)^2 + 4\lambda_1(D_{\alpha})}\bigg).
    $$ Differentiating with respect to $\alpha$
 we get
 
 $$
\frac{d \mu_{2,\alpha}}{d \alpha}|_{\alpha=1}=\frac{d \lambda_{1}(D_\alpha)}{d \alpha}|_{\alpha=1}\frac{4}{N(N-2)}\bigg( \frac{1}{\sqrt{(N-2)^2 + 4\lambda_1(D_{\alpha|\alpha=1})}}+1\bigg).
$$
From Proposition \ref{derneg} we have that 
 $$
 \frac{d \lambda_{1}(D_\alpha)}{d \alpha}|_{\alpha=1}<0,
 $$
from which we get
 $$
 \frac{d}{d\alpha} \mu_{2,\alpha}|_{\alpha=1}\ne0
 $$
 and so the assertion holds.
\end{proof}
Now we prove Theorem \ref{mainth}

\begin{proof}[Proof of Theorem \ref{mainth}]
We apply Theorem \ref{CRT} at the operator $F$ defined in \eqref{Tdef}. It is easy to see that $F(\alpha,U) = 0$ for any $\alpha$. By Lemma \ref{contder}  the operators $\partial_{u}F,\partial_{\alpha}(\partial_{u}F)$ are well defined and continuous from $(0,\alpha^*)\times(K \cap B)$ to $K$. Lemma \ref{ker} says that the kernel of $\partial_u F(1,U)$ is one-dimensional, while Lemma \ref{range} implies that its range has codimension one. Finally, Lemma \ref{transvers} guarantees that the transversality condition (that is condition $(iv)$ of Theorem \ref{CRT}) holds. Therefore, all assumptions of Theorem \ref{CRT} are satisfied. Hence there exists a $C^1$ curve  $(-\epsilon,\epsilon)\ni s\to (\alpha_s,w_s)$, for $\epsilon$ sufficiently small, such that $\alpha_0=1$, $w_0=U$, $\frac{d}{ds}w_s\big|_{s=0}=\psi_2,$ and
$$
\begin{cases}
   \mathcal{L}_{\alpha_s} w_s=w_s^{2^*-1} &\text{in }\Sigma_{D_{1}}\\
    w_s\in \mathcal{{D}}^{1,2}(\Sigma_{D_1})\\
    w_s>0,
    \end{cases}
$$
taking $v_s=\Phi_{1/\alpha}^*w_s$ we conclude.

\end{proof}

Next we prove a global bifurcation result. We use the following Theorem due to Rabinowitz (\cite{K}, Theorem II.3.3) that we state in the case when the parameter interval is bounded.
\begin{theorem}
Assume that $F \in C(X\times J,X), F(x,t) = x + f(x,t)$, where $J=(a,b)$ is a bounded open interval in $\R$, $f :  X\times J\to X$ is completely continuous, and $D_xF(0,\cdot) = I + D_xf(0,\cdot) \in C(J, L(X,X))$. Let $S$ denote the closure of the set of nontrivial solutions of $F(x,t) = 0$ in $X\times J$. Assume that $D_xF(0,t)$ has an odd crossing number at $t= t_0$. Then $(0,t_0) \in S$, and let $C$ be the (connected) component of $S$ to which $(0,t_0)$ belongs. Then:

either
\begin{enumerate}
    \item $C$ is unbounded, 
    
    or
    \item $ C$ contains some $(0,t_1)$, with $t_0 \ne t_1, t_1\in J$

    or

    \item $\bar C\cap (X\times\{a\}\cup X\times\{b\})\ne\emptyset.$
\end{enumerate}
\end{theorem}

Note that the alternative $(3)$ is easily deduced from the proof of Theorem II.3.3 in \cite{K} if the parameter interval is not the whole $\R$.

Applying  this theorem to our case we get:
\begin{theorem}
    Let $S$ denote the closure of the set of nontrivial solutions of $F(\alpha,u) = 0$ in $(0,\alpha^*)\times K$.  Then $(1,U) \in S$, and let $C$ be the (connected) component of $S$ to which $(1,U)$ belongs. Then any $u\in C$ is positive on $\bar\Sigma_D$ and:

either
\begin{enumerate}
    \item $C$ is unbounded, 
    
    or
    \item $ C$ contains some $(\alpha_1,U)$, where $\alpha_1\ne 1$ $\alpha_1\in(0,\alpha^*)$,

    or
    \item $\bar  C\cap (\{0\}\times K\cup\{\alpha^*\}\times K)\ne \emptyset$.
\end{enumerate}
\end{theorem}
\begin{proof}
We have to prove that $f : (0,\alpha^*)\times K\to K$, 

$$f(\alpha,u)=-( \mathcal{L}_\alpha
)^{-1}(u^{2^*-1})
$$
is completely continuous, where $\mathcal L_{\alpha}$ is defined by
$$
\mathcal L_{\alpha}(\Phi_{\alpha}^*v)=\Phi_{\alpha}^*(-\Delta v).
$$
First of all if $u\in K$ then $\mathcal{L}_\alpha^{-1}(u^{2^*-2}):((0,\alpha^*),L^2(u^{2-2^*}))\to \mathcal{D}^{1,2}(\Sigma_{D_1})$ and one gets $f(\alpha,u)\in \mathcal{D}^{1,2}(\Sigma_{D_1})$. Moreover defining  $\tilde f=\Phi^*_{1/\alpha}f$ and $\tilde u=\Phi^*_{1/\alpha}u$ it holds that
$$
\Phi_{\alpha}^*(\Delta\tilde f)=-\mathcal L_{\alpha}(\Phi_{\alpha}^*\tilde f)=u^{2^*-1}=\Phi_{\alpha}^*(\tilde u)
$$
and hence
$$
\tilde f(\alpha,u)=\Delta^{-1}(\tilde u^{2^*-1}).
$$
From the representation formula we get
\beq\label{compineq}
|\tilde f(\alpha,u)|\le C\int_{\Sigma_{D_\alpha}}\frac{1}{|x-y|^{N-2}}\tilde u(y)^{2^*-1}dy\le C\int_{\Sigma_{D_\alpha}}\frac{U^{2^*-1}(y)}{|x-y|^{N-2}}dy=CU(x).
\eeq
This yields
$$
\Phi^*_{1/\alpha}(f)\le|\tilde f(\alpha,u)|\le CU(x)
$$
which implies
\beq\label{bddbubble24}
f\le\Phi^*_{\alpha}( CU(x))=CU(x)
\eeq
and so $f(\alpha,u)\in K$. In particular we have that $f : (0,\alpha^*)\times K\to K$ is  continuous.

We must show that if $(\alpha_m, u_m) \in (0,\alpha^*)\times K $ is bounded, with $ \alpha_m $ converging to some limit $ \alpha $ in $ [0,\alpha^*] $ (up to a subsequence), then $f_m=f(\alpha_m, u_m) $ contains a converging subsequence. 
Since the sequence $(u_m)$ is bounded in the Hilbert space $\mathcal{D}^{1,2}(\Sigma_{D_1})$, up to a subsequence (still denoted $u_m$), we have:
\[
u_m \rightharpoonup u \quad \text{weakly in } \mathcal{D}^{1,2}(\Sigma_{D_1}),
\]
and then $u_m$ converges strongly in $L^2(u^{2^*-2})$. Up to a further subsequence,
\[
u_m(x) \to u(x) \quad \text{a.e. in } \Sigma_{D_1},
\]
then since \( u_m \) is Kelvin invariant for all \( m \in \mathbb{N} \), we conclude that \( u \) is Kelvin invariant as well.
Moreover $\|u_m\|_K$ is bounded which means $|u_m|\le CU$ where $C$ is independent of $m$ and so $|u_m|^{2^*-1}\le C U^{2^*-1}$. Then we have

$$
\int_{\Sigma_{D_1}} \mathcal{L}_{\alpha_m} f_m \varphi-\int_{\Sigma_{D_1}} \mathcal{L}_{\alpha} f \varphi =\int_{\Sigma_{D_1}} u_m^{2^*-1}\varphi-\int_{\Sigma_{D_1}} u^{2^*-1}\varphi\to 0.
$$
From the uniform ellipticity of $\mathcal{L}_{\alpha}$, $\exists C=C(\alpha)$ such that

$$
C\|\nabla (f_m-f)\|\le \int_{\Sigma_{D_1}} \mathcal{L}_{\alpha_m} (f_m-f)(f_m-f) =\int_{\Sigma_{D_1}} (\mathcal{L}_{\alpha_m} (f_m)-\mathcal{L}_{\alpha} (f))(f_m-f) +$$

$$-\int_{\Sigma_{D_1}} (\mathcal{L}_{\alpha_m} (f)-\mathcal{L}_{\alpha} (f))(f_m-f) =\int_{\Sigma_{D_1}} (u_m^{2^*-1}-u^{2^*-1})(f_m-f)-\int_{\Sigma_{D_1}} (\mathcal{L}_{\alpha_m} (f)-\mathcal{L}_{\alpha} (f))(f_m-f)
$$
$$
=I+II.
$$
We consider $I$. We have that $f_m$ satisfies \eqref{bddbubble24}, then $ (u_m^{2^*-1}-u^{2^*-1})(f_m-f)\le C U^{2^*}\in L^1(\Sigma_{D_1})$ and $(u_m^{2^*-1}-u^{2^*-1})(f_m-f)\to 0$ a.e. in $\Sigma_{D_1}$, then from from Lebesgue's dominated convergence theorem we have that
$$
I\to 0.
$$
Moreover, since $\mathcal{L}_\alpha$ is continuous in operator norm and $f_m-f$ is bounded from \eqref{bddbubble24}, we also get
$$
II\to 0.
$$
One concludes that $f(\alpha_m,u_m)\to f(\alpha,u)$ in $\mathcal{D}^{1,2}(\Sigma_{D_1})$. Moreover passing to the limit on $|u_m|\le CU$  yields $u\in K$  and passing to the limit in \eqref{bddbubble24} with $f_m$, yields $f(\alpha,u)\in K$.

Now we have to prove that $f(\alpha_m,u_m)\to f(\alpha,u)$ in $K$, namely that
$$
\sup_{x\in\Sigma_{D_1}}\frac{|f(\alpha_m,u_m)(x)- f(\alpha,u)(x)|}{U(x)}\to 0.
$$
We consider again  $\tilde f=\Phi^*_{1/\alpha}(f)$ and $\tilde u=\Phi^*_{1/\alpha}(u)$ (respectively $\tilde f_m=\Phi^*_{1/\alpha_m}(f)$ and $\tilde u_m=\Phi^*_{1/\alpha_m} (u)$). Then it holds that
$$
\tilde f_m=\Delta^{-1}(\tilde u_m^{2^*-1})
$$
and
$$
\tilde f=\Delta^{-1}(\tilde u^{2^*-1}).
$$
Then from the representation formula applied to $\tilde f$ and $\tilde f_m$ we get 
$$
|\tilde f(\alpha_m,u_m)(x)- \tilde f(\alpha,u)(x)|\le C\int_{ \Sigma_{D_\alpha}}\frac{|\tilde u_m^{2^*-1}(\tilde y)-\tilde u^{2^*-1}(\tilde y)|}{|\tilde x-\tilde y|^{N-2}}d\tilde y\le$$
and from the Hölder inequality
$$\le C\bigg(\int_{\Sigma_{ D_\alpha}}\bigg(\frac{U^{2^*-1-\epsilon
}(\tilde y)}{|\tilde x-\tilde y|^{N-2}}\bigg)^\frac{q}{q-1}dy\bigg)^\frac{q-1}{q}\bigg(\int_{\Sigma_{ D_\alpha}}\bigg(\frac{|\tilde u_m^{2^*-1}-\tilde u^{2^*-1}|U^{\epsilon
}(\tilde y)}{U^{2^*-1}}\bigg)^qd\tilde y\bigg)^\frac{1}{q}
$$
where $\epsilon>0$ will be chosen small and $q>1$ large such that $\epsilon q=2^*$. The ratio on the integral at the right-hand side is bounded by $C^q U^{\epsilon q}= CU^{2^*}\in L^1(\Sigma_{D_\alpha})$, where $C$ is independent of $N$. Lebesgue's dominated convergence theorem then implies that this integral converges to $0$ as $m\to \infty$.
Then by \cite{GGT} Lemma 3.4 and 3.5,  we can prove that
$$
\int_{\Sigma_{ D_\alpha}}\bigg(\frac{U^{2^*-1-\epsilon
}(\tilde y)}{|\tilde x-\tilde y|^{N-2}}\bigg)^\frac{q}{q-1}dy\le\frac{C}{(1+|\tilde x|)^\frac{(N-2)q}{q-1}}=CU^\frac{q}{q-1}(\tilde x),
$$
and we obtain
$$
\sup_{x\in\Sigma_{D_1}}\frac{|\tilde f(\alpha_m,u_m)(x)- \tilde f(\alpha,u)(x)|}{U(x)}\to 0.
$$

from the definition of $\tilde f_m, \tilde f$ we conclude
$$
\sup_{x\in\Sigma_{D_1}}\frac{|f(\alpha_m,u_m)(x)- f(\alpha,u)(x)|}{U(x)}\to 0.
$$
    We are left to prove that $D_uF(\alpha,U)$ has an odd crossing number at $\alpha =1$. Namely we have to prove that $D_uF(\alpha,U)$ is regular for $\alpha\in(1-\delta,1)\cup(1,1 +\delta)$ and the degree $\sigma^<(\alpha):=(-1)^{m_1(\alpha)+...+m_k}(\alpha)$ changes at $\alpha=1$, where $m_1(\alpha),..,m_k(\alpha)$ are the algebraic multiplicities of all negative real eigenvalues of $ \partial_u F(\alpha,U)$ close to $0$.   This follows easily from the fact that $\lambda_1(D_\alpha)=N-1$ is a simple eigenvalue.
    
    Indeed the eigenvalue problem for the linear operator $
    \partial_u F(\alpha,U)=w-(2^*-1)( \mathcal{L}_\alpha)^{-1}U^{2^*-2}w
    $, is:
$$
w-(2^*-1)( \mathcal{L}_\alpha)^{-1}U^{2^*-2}w=\beta w
    $$
    that is equivalent to
    $$
\mathcal{L}_\alpha w=\frac{2^*-1}{1-\beta}U^{2^*-2}w.
    $$
    From Proposition \ref{lemma_autov} and Corollary \ref{cor_autoval} we know that the first eigenvalue is given by
    $$
    \mu_1=\frac{2^*-1}{1-\beta_1}=1\quad\Rightarrow\quad \beta_1=2-2^*<0
    $$
    the second one $\beta_2$ is given by

    $$    
\mu_2=\frac{2^*-1}{1-\beta_2}= \frac{1}{N(N-2)}\sqrt{(N-2)^2 + 4\lambda_1(D_{\alpha})}\bigg( 2 +\sqrt{(N-2)^2 + 4\lambda_1(D_{\alpha})}\bigg)
    $$

    and it holds that
$$
\mu_2=\frac{2^*-1}{1-\beta_2}\to\begin{cases}
<(2^*-1), &\lambda_1(D_{\alpha})<N-1\\
=2^*-1, &\lambda_1(D_{\alpha})=N-1\\
>(2^*-1), &\lambda_1(D_{\alpha})>N-1,
\end{cases}
$$
    namely
      $$    
\beta_2\to\begin{cases}
<0 &\alpha<1\\
=0 &\alpha=1\\
>0 &\alpha>1.
\end{cases}
    $$
    Then if we take $\delta$ small enough, such that  $\mu_3>0$  we get that $D_uF(\alpha,U)$ has an odd crossing number at $\alpha =1$ and we have the three alternatives. The fact that the solutions $u\in C$ are positive follows from the fact that the subset of $C$ made by positive functions on $\bar \Sigma_D$ is closed and open in $C$ and hence coincides with $C$. This is a consequence of the definition of the space $X$ and the fact that the functions are Kelvin invariant (see also Remark \ref{rmkpalla}). 

\end{proof}

\bibliographystyle{plain}

\end{document}